\documentclass[12pt]{article}
\usepackage{amsmath,amsfonts,latexsym,amsthm,amssymb}
\DeclareMathOperator{\const}{const}
\DeclareMathOperator{\dist}{dist}
\DeclareMathOperator{\diag}{diag}

\DeclareMathOperator{\DRL}{D_{0+}^\alpha }
\DeclareMathOperator{\DC}{\mathbb D^{(\alpha )}}
\begin{document}
\newtheorem{teo}{Theorem}
\newtheorem{lem}{Lemma}
\newtheorem{prop}{Proposition}
\newtheorem{cor}{Corollary}
\title{Fractional differential equations: $\alpha$-entire
solutions, regular and irregular singularities}
\author{Anatoly N. Kochubei
\footnote{Partially supported by the Scientific Program of
the National Academy of Sciences of Ukraine, Project No. 0107U002029, by the DFG, Grant 436 UKR 113/87/01, and
by the Ukrainian Foundation
for Fundamental Research, Grant 14.1/003.}
\\ \footnotesize Institute of Mathematics,\\
\footnotesize National Academy of Sciences of Ukraine,\\
\footnotesize Tereshchenkivska 3, Kiev, 01601 Ukraine.\\
\\ \footnotesize E-mail: kochubei@i.com.ua \quad Phone: (380-44)-243-4515 \quad Fax: (380-44)-235-2010}
\date{}
\maketitle
\vspace*{2cm}
Running head:\quad  ``Fractional differential equations''

\vspace*{2cm}
\begin{abstract}
We consider fractional differential equations of order $\alpha \in (0,1)$ for functions of one independent variable $t\in (0,\infty )$ with the Riemann-Liouville and Caputo-Dzhrbashyan fractional derivatives. A precise estimate for the order of growth of $\alpha$-entire solutions is given. An analog of the Frobenius method for systems with regular singularity is developed. For a model example of an equation with a kind of an irregular singularity, a series for a formal solution is shown to be convergent for $t>0$ (if $\alpha$ is an irrational number poorly approximated by rational ones) but divergent in the distribution sense.
\end{abstract}

\vspace{2cm}

\newpage
\section{INTRODUCTION}

Fractional differential equations are widely used for modeling anomalous relaxation and diffusion phenomena; see \cite{EIK,KST} for further references. Meanwhile the mathematical theory of such equations is still in its initial stage. In particular, a systematic development of the analytic theory of fractional differential equations with variable coefficients was initiated only recently, in the paper by
Kilbas, Rivero, Rodr\'\i guez-Germ\'a, and Trujillo \cite{KRRT} (see also Section 7.5 in \cite{KST}). For equations of order $\alpha \in (0,1)$, of the form
\begin{equation}
\left( \DRL u\right) (t)=a(t)u(t),\quad t>0,
\end{equation}
where $\DRL$ is the Riemann-Liouville fractional derivative, or
\begin{equation}
\left( \DC u\right) (t)=a(t)u(t),\quad t>0,
\end{equation}
where $\DC$ is the Caputo-Dzhrbashyan fractional derivative, their main results are as follows. If $a(t)=A(t^\alpha )$, and $A(z)$ is a real function possessing an absolutely convergent Taylor expansion on an interval $|z|<\theta$, then the equation (1) possesses a solution of the form
$$
u(t)=t^{\alpha -1}\sum\limits_{n=0}^\infty a_nt^{\alpha n},\quad 0<t<\theta ,
$$
while the equation (2) has a solution
$$
u(t)=\sum\limits_{n=0}^\infty b_nt^{\alpha n},\quad
0\le t<\theta .
$$
In both cases the solutions are unique, if appropriate initial conditions are prescribed.

Thus, for example, the property of $\alpha$-analyticity of the coefficient $a(t)$ (defined above) implies a similar property of a solution of the equation (2). In fact, we have $u(t)=U(t^\alpha )$ where $U$ is holomorphic in a disk $\{ z\in \mathbb C,\ |z|<\theta \}$. The coefficient $a$ may be complex-valued as well.

The above results open the way for developing a theory of $\alpha$-analytic solutions of fractional differential equations in the spirit of classical analytic theory of ordinary differential equations. Here we give some results in this direction.

If in (2) $a(t)=A(t^\alpha )$ where $A$ is an entire function, then the above results from \cite{KRRT} with $\theta =\infty$ (stated there in a weaker form, only for real arguments of analytic functions, than actually proved) show that the solution of the Cauchy problem for the equation (2) is of the form $u(t)=U(t^\alpha )$ where $U$ is an entire function. Following $\cite{KRRT,KST}$ we call such solutions $\alpha$-entire. In particular, that is true, if $A$ is a polynomial. A natural question is about the order of $U$ (here we investigate this subject just for the equation (2) since its properties are closer to those of ordinary differential equations).
It is known (see, for example, \cite{CGH} or \cite{HKR}) that every nontrivial solution of the equation $u^{(k)}(z)=A(z)u(z)$, $k\in \mathbb N$, with a polynomial coefficient $A$, is an entire function of order $1+\deg (A)/k$. In this paper we prove that the orders of the entire functions $U$ corresponding to solutions of (2) do not exceed $(1+\deg (A))/\alpha$. As $\alpha \to 1$, this agrees with the above differential equation result. On the other hand, if $\deg (A)=0$, that is $A(z)=\lambda$, $\lambda \in \mathbb C$, then $U(z)=E_\alpha (\lambda z)$, where $E_\alpha$ is the Mittag-Leffler function whose order is $1/\alpha$ \cite{EIK,KST}, which shows the exactness of our general estimate.

Next, we investigate systems of fractional equations with regular singularity, that is the equations
\begin{equation}
t^\alpha \left( \DRL u\right) (t)=A(t^\alpha )u(t)
\end{equation}
and
\begin{equation}
t^\alpha \left( \DC u\right) (t)=A(t^\alpha )u(t)
\end{equation}
where $A(z)$ is a holomorphic matrix-function. Under some assumptions, we prove that formal power series solutions of (3) and (4) converge near the origin and develop an analog of the classical Frobenius method of finding a solution. For scalar equations, the latter problem was considered in \cite{KST,RRT}.

Finally, in order to clarify characteristic features of fractional equations with irregular singularity, we study a model example, the equation
\begin{equation}
t^{2\alpha }\left( \DC u\right) (t)=\lambda u(t),\quad \lambda \in \mathbb C,
\end{equation}
where, as before, $0<\alpha <1$. Assuming that $\alpha$ is irrational and satisfies a Diophantine condition (that is $\alpha$ is poorly approximated by rational numbers), we construct a kind of a formal solution of (5) convergent for $t>0$. We prove that the series for the formal solution does not converge in the distribution sense, within a theory of distributions associated with the fractional calculus (see \cite{R,SKM}). Thus, the formal solution $u(t)$ cannot be interpreted as a distribution solution. It is interesting that $u(t)$ is closely connected with a class of analytic functions with irregular behavior introduced by Hardy \cite{H}.

\section{Preliminaries}

{\bf 2.1. Fractional derivatives and integrals \cite{EIK,KST,SKM}.} Let $\alpha \in (0,1)$ be a fixed number. {\it The Riemann-Liouville fractional integral} of order $\alpha$ of a function $\varphi \in L_1(0,T)$ is defined as
$$
\left( I_{0+}^\alpha \varphi \right) (t)=\frac{1}{\Gamma (\alpha )}\int\limits_0^t(t-\tau )^{\alpha -1}\varphi (\tau )\,d\tau ,\quad 0<t\le T.
$$
{\it The Riemann-Liouville fractional derivative} of order $\alpha$ is given by the expression
$$
\left( D_{0+}^\alpha \varphi \right) (t)=\frac{1}{\Gamma (1-\alpha )}\frac{d}{dt}\int\limits_0^t(t-\tau )^{-\alpha }\varphi (\tau )\,d\tau ,\quad 0<t\le T,
$$
that is $\left( \DRL \varphi \right)(t)=\frac{d}{dt}\left( I_{0+}^{1-\alpha }\varphi \right) (t)$, provided the fractional integral $I_{0+}^{1-\alpha }\varphi$ is an absolutely continuous function. If $\varphi$ is defined on the whole half-axis $(0,\infty )$, then $I_{0+}^\alpha \varphi$ and $\DRL \varphi$ are also defined on $(0,\infty )$. Below we will consider just this case.

Note that the Riemann-Liouville derivative is defined for some functions with a singularity at the origin. For example, if $\varphi (t)=t^d$, $d>-1$, then
\begin{equation}
\left( \DRL \varphi \right) (t)=\frac{\Gamma (d+1)}{\Gamma (d+1-\alpha )}t^{d-\alpha },
\end{equation}
so that $\DRL \varphi =0$, if $\varphi (t)=t^{\alpha -1}$. For $\varphi (t)=t^d$, $d>-1$, we have also
\begin{equation}
\left( I_{0+}^\alpha  \varphi \right) (t)=\frac{\Gamma (d+1)}{\Gamma (d+1+\alpha )}t^{d+\alpha }.
\end{equation}

The Riemann-Liouville fractional differentiation and integration are inverse to each other in the following sense. If $\varphi \in L_1(0,T)$, then $\DRL I_{0+}^\alpha  \varphi =\varphi$. The equality $I_{0+}^\alpha  \DRL \varphi =\varphi$ holds under the stronger assumption that $\varphi =I_{0+}^\alpha \psi$ with some $\psi \in L_1(0,T)$. The latter is equivalent to the conditions of absolute continuity of $I_{0+}^{1- \alpha }\varphi$ on $[0,T]$ and the equality $\left( I_{0+}^{1-\alpha }\varphi \right) (0)=0$.

Let a function $\varphi$ be continuous on $[0,T]$ and possess the Riemann-Liouville fractional derivative of order $\alpha$. The function
\begin{equation}
\left( \DC \varphi \right) (t)=\frac{1}{\Gamma (1-\alpha )}\left[ \frac{d}{dt}\int\limits_0^t(t-\tau )^{-\alpha }\varphi (\tau )\,d\tau -t^{-\alpha }\varphi (0)\right] =\left( \DRL \varphi \right) (t)-\frac{1}{t^\alpha \Gamma (1-\alpha )}\varphi (0)
\end{equation}
is called {\it the Caputo-Dzhrbashyan, or regularized, fractional derivative}. If $\varphi$ is absolutely continuous on $[0,T]$, then
\begin{equation}
\left( \DC \varphi \right) (t)=\frac{1}{\Gamma (1-\alpha )}
\int\limits_0^t(t-\tau )^{-\alpha }\varphi' (\tau )\,d\tau .
\end{equation}
In contrast to $\DRL$, $\DC$ is defined only on continuous functions and vanishes on constant functions. In most of physical applications, equations with $\DC$ are used, because a solution of an equation with the Riemann-Liouville derivative typically has a singularity at the origin $t=0$, so that the initial state of a system to be described by the equation is not defined.

Let $v(z)=\sum\limits_{n=0}^\infty c_nz^n$ be an entire function. Consider the Caputo-Dzhrbashyan derivative of the function
$$
\varphi (t)=v(t^\alpha )=\sum\limits_{n=0}^\infty c_nt^{\alpha n}.
$$
It follows from (6) and (8) that
$$
\left( \DC \varphi \right) (t)=\sum\limits_{n=1}^\infty c_n\beta (n)t^{\alpha (n-1)}=\left( \mathfrak D_\alpha v\right) (t^\alpha )
$$
where
\begin{equation}
\beta (n)=\frac{\Gamma (n\alpha +1)}{\Gamma (n\alpha +1-\alpha )},
\end{equation}
and the operator
$$
\left( \mathfrak D_\alpha v\right) (z)=\sum\limits_{n=1}^\infty c_n\beta (n)z^{n-1}
$$
is known as the Gelfond-Leontiev generalized differentiation operator (see \cite{K,SKM}; in fact, $\mathfrak D_\alpha$ is defined for wider classes of functions).

The operator, right inverse to $\mathfrak D_\alpha$, has the form
$$
\left( \mathfrak I_\alpha f\right) (z)=\frac{z}{\Gamma (\alpha )}\int\limits_0^1(1-t)^{\alpha -1}f(zt^\alpha )\,dt
$$
(see Sect. 22.3 in \cite{SKM}). It will be convenient to make a change of variables setting $z=Re^{i\theta }$, $r=Rt^\alpha$. Then
\begin{equation}
\left( \mathfrak I_\alpha f\right) (Re^{i\theta })=\frac{e^{i\theta }}{\alpha \Gamma (\alpha )}\int\limits_0^R
\left[ \left( \frac{R}r\right)^{1/\alpha}-1\right]^{\alpha -1}f(re^{i\theta })\,dr.
\end{equation}

If $f(z)=\sum\limits_{k=0}^\infty f_kz^k$, then
$$
\left( \mathfrak I_\alpha f\right) (z)=\sum\limits_{k=0}^\infty \frac{\Gamma (\alpha k+1)}{\Gamma (\alpha k+1+\alpha )}f_kz^{k+1},
$$
and it is easy to check that
\begin{equation}
\left( \mathfrak I_\alpha \mathfrak D_\alpha v\right) (z)=v(z)-v(0).
\end{equation}

\medskip
{\bf 2.2. A class of distributions.} Spaces of test functions and distributions behaving reasonably under the action of fractional integration operators were introduced by Rubin \cite{R} (for a brief exposition see also \cite{SKM}; both in \cite{SKM} and \cite{R} there are references regarding other approaches and earlier publications in this field). Proceeding from \cite{R}, it is easy to come to a class of distributions where the Caputo-Dzhrbashyan derivative $\DC$ is defined in a natural way.

Let $\mathcal S (0,\infty )$ be the Schwartz space of smooth functions on $[0,\infty )$ with rapid decay at infinity. Denote
$$
\mathcal S_+=\{ \varphi \in \mathcal S (0,\infty ):\ \varphi^{(l)}(0)=0,\ l=0,1,2,\ldots \},
$$
$$
\Phi_+^{1-\alpha }=\left\{ \varphi \in \mathcal S_+:\ \int\limits_0^\infty \varphi (x)x^{1-\alpha -k}\,dx=0,\ k=1,2,\ldots \right\} ,
$$
$$
\Phi_+^{\alpha -1}=\left\{ \varphi \in \mathcal S_+:\ \int\limits_0^\infty \varphi (x)x^{\alpha -1-k}\,dx=0,\ k=0,1,2,\ldots \right\} .
$$
These spaces are interpreted as topological vector spaces with the topologies induced from $\mathcal S (0,\infty )$; see \cite{R} for various descriptions of these topologies including the description by seminorms.

Together with the Riemann-Liouville fractional integration operator $I_{0+}^\alpha$, it is convenient to use the operator
$$
\left( I_-^\alpha \varphi \right) (t)=\frac{1}{\Gamma (\alpha )}\int\limits_t^\infty (\tau -t)^{\alpha -1}\varphi (\tau )\,d\tau ,\quad t>0.
$$
If $\varphi ,\psi$ are sufficiently good functions, for example, if $\varphi \in L_p(0,\infty )$, $\psi \in L_q(0,\infty )$, $p,q>1$, $\dfrac{1}p+\dfrac{1}q=2-\alpha$, then
$$
\int\limits_0^\infty \varphi (x)\left( I_{0+}^{1-\alpha }\psi \right) (x)\,dx=\int\limits_0^\infty \psi (x)\left( I_-^{1-\alpha }\varphi \right) (x)\,dx
$$
(see Section 2.5.1 in \cite{SKM}). We will write this in the notation
$$
\left\langle \varphi ,I_{0+}^{1-\alpha }\psi \right\rangle =\left\langle \psi ,I_-^{1-\alpha }\varphi \right\rangle .
$$

In particular, if $\psi =u'$, where $u\in \mathcal S (0,\infty )$, then
\begin{equation}
\left\langle \varphi ,\DC u\right\rangle =\left\langle u' ,I_-^{1-\alpha }\varphi \right\rangle .
\end{equation}
It is known \cite{R} that $I_-^{1-\alpha }$ acts continuously from $\Phi_+^{1-\alpha }$ onto $\Phi_+^{\alpha -1}$. Therefore the identity (13) can be used to define $\DC u$ as a distribution from $\left( \Phi_+^{1-\alpha }\right)'$, if $u\in C^1(0,\infty )$, and $u'$ has no more than a power-like growth near zero and infinity. This definition agrees with the classical one: if $u$ is continuously differentiable at the origin too, then
\begin{multline*}
\left\langle u' ,I_-^{1-\alpha }\varphi \right\rangle =\frac{1}{\Gamma (1-\alpha )}\int\limits_0^\infty u'(x)\,dx \int\limits_x^\infty (t-x)^{-\alpha }\varphi (t)\,dt \\
=\frac{1}{\Gamma (1-\alpha )}\int\limits_0^\infty \varphi (t)\,dt \int\limits_0^t (t-x)^{-\alpha }u'(x)\,dx=\left\langle \varphi ,\DC u\right\rangle
\end{multline*}
where $\DC u$ is understood in the sense of (9).

A typical example of a function from $\Phi_+^{1-\alpha }$ is the function
$$
\varkappa_\alpha (x)=x^{\alpha -2}\exp \left( -\frac{\log^2x}4\right) \sin \left( \frac{\pi}2\log x\right) .
$$
It is clear that $\varkappa_\alpha \in \mathcal S_+$. Next, if
$$
\varkappa (x)=\exp \left( -\frac{\log^2x}4\right) \sin \left( \frac{\pi}2\log x\right) ,
$$
then we can write explicitly the Mellin transform
$$
\widetilde{\varkappa}(z)=\int\limits_0^\infty x^{z-1}\varkappa (x)\,dx.
$$
Namely, by the formula (4.133.1) from \cite{GR},
\begin{equation}
\widetilde{\varkappa}(z)=2\int\limits_0^\infty e^{-t^2/4}\sinh (zt)\sin \frac{\pi t}2\,dt=2\sqrt{\pi }e^{z^2-\frac{\pi^2}4}\sin \pi z.
\end{equation}
In particular,
$$
\int\limits_0^\infty \varkappa_\alpha (x)x^{1-\alpha -k}\,dx=\widetilde{\varkappa}(-k)=0,\quad k=0,1,2,\ldots ,
$$
so that indeed $\varkappa_\alpha \in \Phi_+^{1-\alpha }$.

In order to have a full concept of a class of distributions, one needs a result regarding density of a space of test functions in some space of integrable functions. This gives a one-to-one correspondence between ordinary functions and distributions they generate. Here we present such a result though it will not be used directly in this paper. For similar properties in other situations see \cite{S}.

\medskip
\begin{prop}
The space $\Phi_+^{1-\alpha }$ is dense in $L_p((0,\infty ),t^{-1}dt)$, $1\le p<\infty$.
\end{prop}

\medskip
{\it Proof}. Let $R=\int\limits_0^\infty \varkappa_\alpha (t)t^{-1}dt$. We have
$$
R=\widetilde{\varkappa }(\alpha -2)=2\sqrt{\pi }\exp \left( (\alpha -2)^2-\frac{\pi^2}4\right) \sin \pi \alpha >0.
$$
Denote
$$
z_N(x)=\frac{N}R\varkappa_\alpha (x^N),\quad N=1,2,\ldots .
$$
Then
\begin{equation}
\int\limits_0^\infty z_N(x)x^{-1}dx=\frac{1}R\int\limits_0^\infty \varkappa_\alpha (t)t^{-1}dt=1.
\end{equation}

Suppose that $f\in L_p((0,\infty ),t^{-1}dt)$. Consider the so-called Mellin convolution
$$
\left( z_N*_Mf\right) (t)=\int\limits_0^\infty z_N(\tau )f\left( \frac{t}\tau \right) \tau^{-1}d\tau ,
$$
that is actually the convolution on the multiplicative group $(0,\infty )$ (note that $t^{-1}dt$ is a Haar measure on that group). Obviously, $z_N*_Mf\in \Phi_+^{1-\alpha }$. Denote by $\| \cdot \|_p$ the norm in $L_p((0,\infty ),\linebreak t^{-1}dt)$.

Using (15) we can write
$$
\left( z_N*_Mf\right) (t)-f(t)=\int\limits_0^\infty z_N(\tau )\left[ f\left( \frac{t}\tau \right) -f(t)\right] \tau^{-1}d\tau .
$$
By the generalized Minkowski inequality,
$$
\left\| z_N*_Mf-f\right\|_p\le \int\limits_0^\infty |z_N(\tau )|\left\{ \int\limits_0^\infty \left| f\left( \frac{t}\tau \right) -f(t)\right|^pt^{-1}dt\right\}^{1/p}d\tau .
$$

Next we use the $L_p$-continuity of shifts on the multiplicative group $(0,\infty )$ (see \cite{HR}, Theorem 20.4). For any $\varepsilon >0$, there exists $\delta >0$, such that
$$
\left\{ \int\limits_0^\infty \left| f\left( \frac{t}\tau \right) -f(t)\right|^pt^{-1}dt\right\}^{1/p}<\varepsilon ,
$$
if $|\tau -1|<\delta$. Thus,
\begin{equation}
\left\| z_N*_Mf-f\right\|_p\le 2\| f\|_p\int\limits_{|\tau -1|\ge \delta }|z_N(\tau )|\tau^{-1}d\tau +\varepsilon \int\limits_{|\tau -1|<\delta }|z_N(\tau )|\tau^{-1}d\tau .
\end{equation}

Note that
$$
\int\limits_0^\infty |z_N(\tau )|\tau^{-1}d\tau =\frac{N}R
\int\limits_0^\infty \left| \varkappa_\alpha (x^N)\right| x^{-1}dx=\frac{1}R\int\limits_0^\infty \left| \varkappa_\alpha (t)\right| t^{-1}dt=C_1
$$
where the constant $C_1$ does not depend on $N$. On the other hand,
$$
|z_N(x)|\le C_2Nx^{N(\alpha -2)}\exp \left( -\frac{N^2\log^2x}4\right) .
$$
If $|x-1|\ge \delta$, then $\log^2x\ge b>0$, so that
$$
\int\limits_{1+\delta }^\infty |z_N(x)|x^{-1}dx\le C_3e^{-\frac{N^2b}4}(1+\delta )^{N(\alpha -2)}\to 0,
$$
as $N\to \infty$, and (for $\delta <1$)
\begin{multline*}
\int\limits_0^{1-\delta }|z_N(x)|x^{-1}dx\le C_2N\int\limits_{-\infty }^{\log (1-\delta )}\exp \left\{ Nt(\alpha -2)-\frac{N^2t^2}4\right\} \,dt \\
=C_2\int\limits_{-\infty }^{N\log (1-\delta )}\exp \left\{ s(\alpha -2)-\frac{s^2}4\right\} \,ds\to 0,
\end{multline*}
as $N\to \infty$.

As a result, we see that, if $N$ is large enough, the first summand in (16) does not exceed $2\|f\|_p\varepsilon$, while the second $\le C_1\varepsilon$. Thus,
$$
\left\| z_N*_Mf-f\right\|_p\longrightarrow 0,\quad \text{as $N\to \infty$.} \qquad \blacksquare
$$

\bigskip
{\bf 2.3. On ratios of the Gamma functions.} We will often use the function
\begin{equation}
\rho (t)=\frac{\Gamma (t+1)}{\Gamma (t+1-\alpha )},\quad -1<t<\infty .
\end{equation}
Here we collect some of its properties.

If $t>\alpha -1$, the integral representation
\begin{equation}
\frac{1}{\rho (t)}=\frac{1}{\Gamma (\alpha )}\int\limits_0^\infty e^{-st}e^{(\alpha -1)s}\left( 1-e^{-s}\right)^{\alpha -1}ds
\end{equation}
holds (see Chapter 4 in \cite{Ol}). It follows from (18) that the function $t\mapsto \dfrac{1}{\rho (t)}$ is strictly monotone decreasing and $\dfrac{1}{\rho (t)}\to 0$, as $t\to \infty$. Since $\Gamma (t+1-\alpha )$ has a pole at $t=\alpha -1$, it is seen from (17) that $\rho (t)\to 0$, as $t\to \alpha -1+0$.

If $-1<t<\alpha -1$, then (by a well-known identity for the Gamma function)
$$
\rho (t)=\frac{\Gamma (-t+\alpha )}{\Gamma (-t)}\cdot \frac{\sin \pi (t+1)}{\sin \pi (t+1-\alpha )}.
$$
The integral representation for the ratio of the Gamma functions \cite{Ol} leads, after an elementary investigation, to the conclusion that $\rho (t)$ is strictly monotone increasing from $-\infty$ to 0.

Thus, we conclude that on the interval $(-1,\infty )$ the function $\rho (t)$ is strictly monotone increasing from $-\infty$ to $\infty$. The inverse function $\gamma (\lambda )$ solving the equation
$$
\rho (t)=\lambda ,\quad \lambda \in \mathbb R,
$$
is a well-defined continuous function. Note that $\rho (0)=\dfrac{1}{\Gamma (1-\alpha )}$, so that $\gamma (\lambda ) \ge 0$, if $\lambda \ge \dfrac{1}{\Gamma (1-\alpha )}$.

It is known (\cite{Ol}, Chapter 4) that
\begin{equation}
\frac{\Gamma (t+a)}{\Gamma (t+b)}\sim t^{a-b}(1+O(t^{-1})),\quad t\to \infty ,
\end{equation}
if $b>a$. In particular,
$$
\rho (t)\sim t^\alpha , \quad \text{as $t\to \infty$}.
$$

For the sequence $\beta (n)$ defined in (10) and appearing in the definition of the Gelfond-Leontiev generalized differentiation operator $\mathfrak D_\alpha$, we have $\beta (n)=\rho (\alpha n)$, so that the above asymptotics implies the relation
\begin{equation}
\beta (n)\sim Cn^\alpha ,\quad n=0,1,2,\ldots .
\end{equation}

\section{$\alpha$-Entire Solutions}

Let us consider $\alpha$-entire solutions of the equation (2) with $a(t)=A(t^\alpha )$ where $A$ is a polynomial of degree $m\ge 0$. We assume the initial condition $u(0)=u_0$.

\medskip
\begin{teo}
Under the above assumptions, the solution $u(t)$ of the equation (2) has the form $u(t)=v(t^\alpha )$ where $v$ is an entire function whose order does not exceed $(1+m)/\alpha$.
\end{teo}

\medskip
{\it Proof}. Seeking the function $v$,
we have $\left( \mathfrak D_\alpha v\right) (z)=a(z)v(z)$. Let us apply the operator $\mathfrak I_\alpha$ (see (11)) to both sides of this equality. We get from (11) and (12) that
$$
v(Re^{i\theta })-v(0)=\frac{e^{i\theta }}{\alpha \Gamma (\alpha )}\int\limits_0^R\left[ \left( \frac{R}r\right)^{1/\alpha}-1\right]^{\alpha -1}a(re^{i\theta })v(re^{i\theta })\,dr,
$$
which implies the inequality
$$
|v(Re^{i\theta })|\le |v(0)|+C\int\limits_0^R\left[ \left( \frac{R}r\right)^{1/\alpha}-1\right]^{\alpha -1}|a(re^{i\theta })||v(re^{i\theta })|\,dr
$$
(here and below we denote by the same letter $C$ various positive constants).

We have the asymptotic relations
$$
t^{\frac{1}\alpha}-1\sim \frac{1}\alpha (t-1),\quad \text{as $t\to 1+0$};
$$
$$
\left( t^{\frac{1}\alpha}-1\right)^{\alpha -1}\sim t^{\frac{\alpha -1}\alpha},\quad \text{as $t\to \infty$}.
$$
Therefore
$$
\left( t^{\frac{1}\alpha}-1\right)^{\alpha -1}\le C(t-1)^{\alpha -1}t^{2-\frac{1}\alpha -\alpha },\quad t\ge 1,
$$
so that
$$
|v(Re^{i\theta })|\le |v(0)|+CR^{2-\frac{1}\alpha -\alpha }\int\limits_0^R(R-r)^{\alpha -1}r^{\frac{1}\alpha -1}|a(re^{i\theta })||v(re^{i\theta })|\,dr.
$$

Since $a(z)$ is a polynomial of degree $m$, we get
$$
|v(Re^{i\theta })|\le |v(0)|+CR^{2-\frac{1}\alpha -\alpha }\int\limits_0^R(R-r)^{\alpha -1}r^{\frac{1}\alpha -1+m}|v(re^{i\theta })|\,dr
$$
where $C$ does not depend on $R,\theta$. Fixing $\theta$ and denoting
$$
w(r)=\frac{|v(re^{i\theta })|}{r^{2-\frac{1}\alpha -\alpha }}
$$
we come to the inequality
\begin{equation}
w(R)\le \frac{|v(0)|}{R^{2-\frac{1}\alpha -\alpha }}+C\int\limits_0^R(R-r)^{\alpha -1}r^{1-\alpha +m}w(r)\,dr.
\end{equation}
Now we are in a position to apply Henry's theorem (see Lemma 7.1.2 from \cite{He}), which states that the inequality (21) implies the inequality
$$
w(R)\le \frac{|v(0)|}{R^{2-\frac{1}\alpha -\alpha }}\mathcal E_{\alpha ,2-\alpha +m}(CR)
$$
where $\mathcal E_{\alpha ,\sigma}(s)$ is a certain function admitting the estimate
$$
\mathcal E_{\alpha ,\sigma}(s)\le Cs^{\frac{1}2\left( \frac{\alpha +\sigma -1}\alpha -\sigma \right)} \exp \left( \frac{\alpha}{\alpha +\sigma -1}s^{\frac{\alpha +\sigma -1}\alpha }\right) .
$$
Thus,
$$
|v(Re^{i\theta })|\le C\exp \left( \mu R^{\frac{1+m}\alpha }\right)
$$
for some $\mu \ge 0$, as desired. $\qquad \blacksquare$

\section{Regular Singularity}

{\bf 4.1. Formal and $\alpha$-analytic solutions.} Let us consider systems of equations of the form
\begin{equation}
t^\alpha \left( \DRL u\right) (t)=A(t^\alpha )u(t)
\end{equation}
where
$$
A(z)=A_0+\sum\limits_{m=1}^\infty A_mz^m,
$$
$A_m$ are $n\times n$ complex matrices and
$$
\| A_m\|\le M\mu^m\ (\mu >0),\quad m=0,1,2,\ldots .
$$

Suppose we have a formal series
\begin{equation}
u(t)=\sum\limits_{k=0}^\infty u_kt^{\alpha k},\quad u_k\in \mathbb C^n.
\end{equation}
Let us substitute the series (23) formally into (22). We get, in accordance with (6), that
$$
\sum\limits_{k=0}^\infty \beta (k)u_kt^{\alpha k}=\sum\limits_{m,k=0}^\infty A_mu_kt^{\alpha (m+k)}
$$
where $\beta (k)$ is the sequence (10). Collecting and comparing the terms we find that
$$
\beta (l)u_l=\sum\limits_{k=0}^l A_ku_{l-k},\quad l=0,1,2,\ldots ,
$$
or, equivalently,
\begin{equation}
A_0u_0=\frac{1}{\Gamma(1-\alpha )}u_0;
\end{equation}
\begin{equation}
\left[ A_0-\beta (l)\right] u_l=-\sum\limits_{k=1}^lA_ku_{l-k},\quad l\ge 1.
\end{equation}

It is natural to call the formal series (23) {\it a formal solution} of the system (22) if the relations (24), (25) hold.

\medskip
\begin{prop}
If a formal series (23) is a formal solution of the system (22), then the series (23) is absolutely convergent on some neighbourhood of the origin.
\end{prop}

\medskip
{\it Proof}. It follows from (19) that
$$
\| \left[ A_0-\beta (l)\right]^{-1}\|\le Cl^{-\alpha },\quad l\ge l_0.
$$
In particular, we may assume that
$$
\| \left[ A_0-\beta (l)\right]^{-1}\|\le 1,\quad l\ge l_0.
$$
Considering, if necessary, $\lambda u(t)$ instead of $u(t)$, with $|\lambda |$ small enough, we may assume that $\|u_0\|\le 1$.

Let us choose so big $r>0$ that $\|u_l\|\le r^l$ for $l\le l_0$ and
$$
M\sum\limits_{k=1}^\infty \left( \frac{\mu }r\right)^k\le 1.
$$
Then
$$
\|u_l\|\le r^l\quad \text{for all $l$.}
$$

Indeed, if this inequality is proved up to some value of $l\ge l_0$, then
$$
\|u_{l+1}\|\le \left\| \sum\limits_{k=1}^{l+1}A_ku_{l+1-k}\right\|\le M\sum\limits_{k=1}^{l+1}\mu^kr^{l+1-k}=Mr^{l+1}\sum\limits_{k=1}^{l+1} \left( \frac{\mu }r\right)^k\le r^{l+1},
$$
and the above inequality implying local convergence in (23) has been proved. $\qquad \blacksquare$

\medskip
The above arguments remain valid for systems of the form
\begin{equation}
t^\alpha \left( \DC u\right) (t)=A(t^\alpha )u(t).
\end{equation}
The only difference is that, instead of (24), we get the relation $A_0u_0=0$, just as in the classical case (see \cite{Hart}).

\bigskip
{\bf 4.2. Model scalar equations.} Consider the equation
\begin{equation}
t^\alpha \left( \DRL \varphi \right) (t)=\lambda \varphi (t),\quad \lambda \in \mathbb R.
\end{equation}
By the relation (6), a solution of the equation (27) is $\const \cdot t^{\gamma (\lambda )}$ where $\gamma$ (the inverse function to $\rho$) was defined in Section 2.3. For example, if $\lambda =0$, then we have $\gamma (0)=\alpha -1$.

If we consider an equation similar to (27), but with the Caputo-Dzhrbashyan derivative, that is
\begin{equation}
t^\alpha \left( \DC \varphi \right) (t)=\lambda \varphi (t),
\end{equation}
then the constant function is a solution of (28) for $\lambda =0$. Suppose that $\lambda \ne 0$, and $\varphi$ is a solution of (28), that is $\varphi \in C[0,T]$, the function
$$
\left( I_{0+}^{1-\alpha }\varphi \right) (t)=\frac{1}{\Gamma (1-\alpha )}\int\limits_0^t(t-\tau )^{-\alpha }\varphi (\tau )\,d\tau
$$
is absolutely continuous, and (28) is satisfied with
$$
\left( \DC \varphi \right) (t)=\frac{d}{dt}\left( I_{0+}^{1-\alpha }\varphi \right) (t)-\dfrac{1}{t^\alpha \Gamma (1-\alpha )}\varphi (0).
$$
We have
$$
\left( I_{0+}^{1-\alpha }\varphi \right) (t)=\frac{t^{1-\alpha }}{\Gamma (1-\alpha )}\int\limits_0^1(1-s)^{-\alpha }\varphi (st)\,ds,
$$
and since $\varphi$ is continuous,
$$
\left( I_{0+}^{1-\alpha }\varphi \right) (t)\longrightarrow 0,\quad \text{as $t\to +0$.}
$$

It is known \cite{SKM} that in these circumstances $I_{0+}^\alpha \DRL \varphi =\varphi$. Note also that $I_{0+}^\alpha$ transforms the function $t^{-\alpha }$ into the constant $\Gamma (1-\alpha )$. Dividing the equation (28) by $t^\alpha$ and applying $I_{0+}^\alpha$ to both sides, we find that
$$
\varphi (t)-\varphi (0)=\frac{\lambda }{\Gamma (\alpha )}\int\limits_0^t(t-\tau )^{-1+\alpha }\tau^{-\alpha }\varphi (\tau )\,d\tau =\frac{\lambda }{\Gamma (\alpha )}\int\limits_0^1(1-s)^{-1+\alpha }s^{-\alpha }\varphi (ts)\,ds.
$$
Passing to the limit, as $t\to 0$, and taking into account the continuity of $\varphi$, we obtain the identity
$$
\frac{\lambda \varphi (0)}{\Gamma (\alpha )}\int\limits_0^1(1-s)^{-1+\alpha }s^{-\alpha }ds=0,
$$
whence $\varphi (0)=0$.

Thus, for $\lambda \ne 0$, the equation (28) is equivalent to (27), if (27) is considered for continuous functions vanishing at the origin. The power solution $Ct^{\gamma (\lambda )}$ belongs to this class, if $\lambda >\frac{1}{\Gamma (1-\alpha )}$. It may be instructive to see these solutions, satisfying the equation
$$
\varphi (t)=\frac{\lambda }{\Gamma (\alpha )}\int\limits_0^t(t-\tau )^{-1+\alpha }\tau^{-\alpha }\varphi (\tau )\,d\tau ,
$$
as examples of non-uniqueness of solutions of linear Volterra integral equations occurring due to the singularity of a kernel.

\bigskip
{\bf 4.3. Systems with good spectrum.} Let us consider the equation (3) with $A(z)=\sum\limits_{m=0}^\infty A_mz^m$ where $A_m$ are complex $n\times n$ matrices, the matrix $A_0$ is Hermitian, and the series converges on a neighbourhood of the origin. Without restricting generality, we may assume that
$$
A_0=\diag (\lambda_1,\ldots ,\lambda_n),\quad \lambda_1,\ldots ,\lambda_n\in \mathbb R.
$$

Following the classical method (see, for example, \cite{B}) we look for a matrix-valued solution ({\it a fundamental solution}) of the equation (3), in the form
\begin{equation}
u(t)=S(t^\alpha )\psi (t)
\end{equation}
where $\psi (t)=\diag \left( t^{\gamma (\lambda_1)},\ldots , t^{\gamma (\lambda_n)}\right)$, $S(z)=\sum\limits_{\nu =0}^\infty \sigma_\nu z^\nu$, $\sigma_\nu$ ($\nu \ge 1$) are some unknown matrices, $\sigma_0=I$.

We have
$$
u(t)=\sum\limits_{\nu =0}^\infty \sigma_\nu \diag \left( t^{\gamma (\lambda_1)+\alpha \nu},\ldots ,t^{\gamma (\lambda_n)+\alpha \nu}\right),
$$
whence
\begin{equation}
t^\alpha \left( \DRL u\right) (t)=\sum\limits_{\nu =0}^\infty \sigma_\nu R_\nu \diag \left( t^{\gamma (\lambda_1)+\alpha \nu},\ldots ,t^{\gamma (\lambda_n)+\alpha \nu}\right)
\end{equation}
where
$$
R_\nu u =\diag \left( \frac{\Gamma (\gamma (\lambda_1)+\alpha \nu +1)}{\Gamma (\gamma (\lambda_1)+\alpha \nu +1-\alpha )},\ldots ,\frac{\Gamma (\gamma (\lambda_n)+\alpha \nu +1)}{\Gamma (\gamma (\lambda_n)+\alpha \nu +1-\alpha )}\right) .
$$
On the other hand,
\begin{equation}
A(t^\alpha )u(t)=\sum\limits_{\nu =0}^\infty \left( \sum\limits_{m=0}^\nu A_m\sigma_{\nu -m}\right)\diag \left( t^{\gamma (\lambda_1)+\alpha \nu},\ldots ,t^{\gamma (\lambda_n)+\alpha \nu}\right).
\end{equation}

Note that $R_0=\diag (\lambda_1,\ldots ,\lambda_n)=A_0$, and since $\sigma_0=I$, the coefficients corresponding to $\nu =0$ in (30) and (31) coincide. Comparing the rest of the coefficients, we obtain the following system of equations for the matrices $\sigma_k$:
\begin{equation}
\sigma_kR_k-A_0\sigma_k=\sum\limits_{l=0}^{k-1}A_{k-l}\sigma_l,\quad k\ge 1.
\end{equation}

For each $k$, the matrix equation (32) for $\sigma_k$ has a unique solution if the spectra of the matrices $R_k$ and $A_0$ are disjoint (see Appendix A.1 in \cite{B}), that is
\begin{equation}
\frac{\Gamma (\gamma (\lambda_i)+\alpha k+1)}{\Gamma (\gamma (\lambda_i)+\alpha k+1-\alpha )}\ne \lambda_j
\end{equation}
for all $i,j\in \{ 1,\ldots ,n\}$, or, equivalently, since the left-hand side of (33) equals $\rho (\gamma (\lambda_i)+\alpha k)$,
$$
\gamma (\lambda_j)-\gamma (\lambda_i)\ne \alpha k,\quad \text{for all $i,j\in \{ 1,\ldots ,n\}$}.
$$

We call our system (3) {\it a system with good spectrum}, if
\begin{equation}
\gamma (\lambda_j)-\gamma (\lambda_i)\notin \alpha \mathbb N, \quad \text{for all $i,j\in \{ 1,\ldots ,n\}$}.
\end{equation}
This definition extends the classical one \cite{B}, since for $\alpha =1$ we would have $\rho (t)=\gamma (t)=t$, and the condition (34) would mean that the eigenvalues of $A_0$ must not differ by a natural number.

\medskip
\begin{teo}
If a system (3) has a good spectrum, then it possesses a fundamental solution (29) where the series for $S(z)$ has a positive radius of convergence.
\end{teo}

\medskip
{\it Proof}. By the asymptotic relation (19),
$$
\frac{\Gamma (\gamma (\lambda_j)+\alpha k+1)}{\Gamma (\gamma (\lambda_j)+\alpha k+1-\alpha )}\sim (\alpha k)^\alpha \left( 1+O(k^{-1})\right),\quad k\to \infty ,
$$
for all $j=1,\ldots ,n$. Therefore
\begin{equation}
(\alpha k)^{-\alpha }R_k=I+O(k^{-1}), \quad k\to \infty .
\end{equation}

Let us divide both sides of the equation (32) by $(\alpha k)^\alpha$. The resulting equation, considered as a system of scalar equations for $n^2$ elements of the matrix $\sigma_k$, has the coefficients bounded in $k$ and the determinant, which is different from zero for each $k$ and tends to 1, as $k\to \infty$. This implies the estimate
\begin{equation}
\|\sigma_k\| \le ak^{-\alpha }\left\| \sum\limits_{l=0}^{k-1}A_{k-l}\sigma_l\right\| ,\quad k\ge 1,
\end{equation}
where the constant $a>0$ does not depend on $k$. It follows from (35), (36), and the convergence near the origin of the power series for $A(z)$ that
$$
\|\sigma_k\| \le a_1k^{-\alpha }\sum\limits_{l=0}^{k-1}b^{k-l}\|\sigma_l\| ,\quad k\ge 1,
$$
where $a_1$ and $b$ are positive constants independent of $k$.

Define a sequence $\{ s_k\}_0^\infty$ of positive numbers, setting $s_0=1$,
$$
s_k=a_1k^{-\alpha }\sum\limits_{l=0}^{k-1}b^{k-l}s_l, \quad k\ge 1.
$$
The induction on $k$ yields the inequality $\| \sigma_k\| \le s_k$ for all $k\ge 0$. On the other hand,
\begin{multline*}
s_{k+1}=a_1(k+1)^{-\alpha }\sum\limits_{l=0}^kb^{k+1-l}s_l
=\frac{(k+1)^{-\alpha}}{k^{-\alpha }}\left[ a_1k^{-\alpha }\left( b\sum\limits_{l=0}^{k-1}b^{k-l}s_l+bs_k\right) \right] \\
=\frac{(k+1)^{-\alpha}}{k^{-\alpha }}\left( bs_k+a_1k^{-\alpha }bs_k\right) =\frac{(k+1)^{-\alpha}b}{k^{-\alpha }}\left( 1+a_1k^{-\alpha}\right) s_k.
\end{multline*}
Therefore
\begin{equation}
\frac{s_k}{s_{k+1}}\longrightarrow b^{-1},\quad \text{as $k\to \infty$.}
\end{equation}

It follows from (37) that the series $\sum\limits_{k=0}^\infty s_kz^k$ has the convergence radius $b^{-1}$ (see Section 2.6 in \cite{WW}). Moreover, the series $\sum\limits_{\nu =0}^\infty \sigma_\nu z^\nu$ converges for $|z|<b^{-1}$. $\qquad \blacksquare$

\medskip
For the equation (26), a similar construction is valid, if we assume that $\lambda_1,\ldots ,\lambda_n\ge \dfrac{1}{\Gamma (1-\alpha )}$.

\section{Irregular Singularity: an Example}

{\bf 5.1. A formal solution.} In this section we construct a solution, in a sense to be specified, of the equation (5). Looking at classical first order equations, corresponding formally to $\alpha =1$, we have to consider the equation $t^2y'(t)=\lambda y(t)$ whose solution is $y(t)=\exp (-\lambda t^{-1})$. Therefore it is natural to seek a solution of the equation (5) in the form
\begin{equation}
u(t)=\sum\limits_{n=0}^\infty c_nt^{-n\alpha},\quad c_n\in \mathbb C.
\end{equation}

A fractional derivative of any term in (38) with $n>\alpha^{-1}$ does not make sense classically. However we may apply the distribution theory from Sect. 2.2. Below we understand the fractional derivative $\DC$ in the sense of (13).

\medskip
\begin{prop}
$\mathrm{(i)}$ If $\mu <0$, $\mu \ne -1,-2,\ldots$, then
\begin{equation}
\DC t^\mu =\frac{\Gamma (\mu +1)}{\Gamma (\mu +1-\alpha )}t^{\mu -\alpha }.
\end{equation}

$\mathrm{(ii)}$ If $k$ is a natural number, then
\begin{equation}
\DC t^{-k}=\frac{(-1)^{k-1}}{(k-1)!\Gamma (-k+1-\alpha )}t^{-k -\alpha }\log t.
\end{equation}
\end{prop}

\medskip
{\it Proof}. Let $\varphi \in \Phi_+^{1-\alpha }$. By (13),
\begin{equation}
\left\langle \DC t^\mu ,\varphi (t)\right\rangle =\mu \left\langle t^{\mu -1},\left( I_-^{1-\alpha} \varphi \right) (t)\right\rangle .
\end{equation}
It is clear that the right-hand side of (41) is an entire function of $\mu$. For $\mu >0$, by virtue of (7),
$$
\left\langle t^{\mu -1},\left( I_-^{1-\alpha} \varphi \right) (t)\right\rangle =\left\langle I_{0+}^{1-\alpha}t^{\mu -1}, \varphi (t)\right\rangle =\frac{\Gamma (\mu )}{\Gamma (\mu +1-\alpha )}\left\langle t^{\mu -\alpha },\varphi (t)\right\rangle .
$$
For $\mu <0$, $\mu \ne -1,-2,\ldots$, the analytic continuation gives the equality (39).

Next, consider the entire function
$$
F(\mu )=\langle t^{\mu -\alpha },\varphi (t)\rangle ,\quad \mu \in \mathbb C.
$$
Note that $F(-k)=0$, $k\in \mathbb N$, by the definition of the space $\Phi_+^{1-\alpha }$. We have
$$
F'(\mu )=\left\langle \log t\cdot e^{(\mu -\alpha )\log t},\varphi \right\rangle .
$$
In particular,
$$
F'(-k)=\int\limits_0^\infty t^{-\alpha -k}\log t\cdot \varphi (t)\,dt.
$$

As $\mu$ belongs to a small neighbourhood of the point $-k$, $F(\mu )=F'(-k)(\mu +k)+o(\mu +k)$. Since the residue of $\Gamma (\mu +1)$ at $\mu =-k$ equals $\dfrac{(-1)^{k-1}}{(k-1)!}$ (see Section 4.4.1 of \cite{T}), we see that the function $\Gamma (\mu +1)F(\mu )$ is holomorphic at $\mu =-k$ (in fact, it is entire) and
$$
\Gamma (\mu +1)F(\mu )|_{\mu =-k}=\frac{(-1)^{k-1}}{(k-1)!}
\int\limits_0^\infty t^{-\alpha -k}\log t\cdot \varphi (t)\,dt.
$$
Now the equality (39) implies (40). $\qquad \blacksquare$

\medskip
Returning to (38), we will formally apply $\DC$ termwise and find the coefficients $c_n$ comparing resulting terms in (5). It is clear from the equality (40) that such a procedure would fail if some of the numbers $n\alpha$ are integers. Thus, we have to assume that $\alpha$ {\it is irrational}. Using (39) we find that (formally)
$$
t^{2\alpha }\left( \DC u\right) (t)=\sum\limits_{n=1}^\infty c_n\frac{\Gamma (-n\alpha +1)}{\Gamma (-n\alpha +1-\alpha)}t^{-n\alpha +\alpha}.
$$
Substituting this into (5) we come to the recurrence relation
$$
c_{n+1}=\lambda \frac{\Gamma (1-(n+2)\alpha )}{\Gamma (1-(n+1)\alpha )}c_n,\quad n\ge 0,
$$
and it is easy to find by induction that
$$
c_n=\lambda^n\frac{\Gamma (1-(n+1)\alpha )}{\Gamma (1-\alpha )}c_0.
$$

Thus, we have found a formal solution
\begin{equation}
u(t)=\frac{c_0}{\Gamma (1-\alpha )}\sum\limits_{n=0}^\infty \Gamma (1-(n+1)\alpha )t^{-n\alpha }
\end{equation}
of the equation (5). Using the identity
$$
\Gamma (z)\Gamma (1-z)=\frac{\pi }{\sin \pi z}
$$
we can rewrite (42) in the form
\begin{equation}
u(t)=\frac{c_0\pi }{\Gamma (1-\alpha )}\sum\limits_{n=0}^\infty \frac{1}{\sin (\pi (n+1)\alpha )}\frac{t^{-n\alpha }}{\Gamma ((n+1)\alpha )}.
\end{equation}

\bigskip
{\bf 5.2. The convergence problem.} The convergence of the series (43) depends on the arithmetic properties of the irrational number $\alpha$. It was shown by Hardy \cite{H} that $\alpha$ can be chosen in such a way (to be well approximated by rational numbers) that the series (43) would diverge for small values of $t$.

An irrational number $\alpha \in (0,1)$ is said to be {\it poorly approximated by rational numbers} if there exist such $\varepsilon >0$, $c>0$ that for any rational number $\dfrac{p}q$, $p,q\in \mathbb N$,
\begin{equation}
\left| \alpha -\frac{p}q\right| \ge cq^{-2-\varepsilon }.
\end{equation}
By the Thue-Siegel-Roth theorem (see \cite{Sch}) such are all algebraic numbers.

The first statement of the next theorem is actually contained already in the paper \cite{H}.

\medskip
\begin{teo}
If $\alpha$ is poorly approximated by rational numbers, then the series (43) converges for any $t>0$. However this series diverges in the space of distributions $\left( \Phi_+^{1-\alpha }\right)'$.
\end{teo}

\medskip
{\it Proof}. It follows from (44) that
$$
|q\alpha -p|\ge cq^{-1-\varepsilon }\quad \text{for any $p\in \mathbb N$,}
$$
so that
$$
\dist ((n+1)\alpha ,\mathbb Z_+)\ge c(n+1)^{-1-\varepsilon },
$$
and taking $l\in \mathbb Z_+$, such that $|(n+1)\alpha -l|\le \dfrac12$, we find that
$$
|\sin (\pi (n+1)\alpha )|=|\sin (\pi ((n+1)\alpha -l)|\ge 2|(n+1)\alpha -l|\ge 2\dist ((n+1)\alpha ,\mathbb Z_+)\ge 2c(n+1)^{-1-\varepsilon }.
$$
Using the Stirling formula we obtain that the series (43) converges for each $t>0$.

To prove the second assertion, consider $\left\langle t^{-\alpha n},\varkappa_\alpha (t)\right\rangle$ where the function $\varkappa_\alpha$ was defined in Sect. 2.2. We have, by (14), that
\begin{multline*}
\left\langle t^{-\alpha n},\varkappa_\alpha (t)\right\rangle =\int\limits_0^\infty t^{-\alpha n}\varkappa_\alpha (t)\,dt=\widetilde{\varkappa_\alpha}(1-\alpha n)=\widetilde{\varkappa}(\alpha -1-\alpha n)\\
=2\sqrt{\pi }\exp \left( (\alpha -1-\alpha n)^2-\frac{\pi^2}4\right) \sin (\pi (\alpha -1-\alpha n)).
\end{multline*}

Now we can give a lower estimate of the coefficients in the series (43) understood in the distribution sense: there are such $a,a_1,b>0$ that
\begin{multline*}
\left| \frac{1}{\sin (\pi (n+1)\alpha )}\cdot \frac{\left\langle t^{-n\alpha },\varkappa_\alpha (t)\right\rangle }{\Gamma ((n+1)\alpha )}\right| \ge a(n+1)^{-b(n+1)}|\sin (\pi \alpha (n-1))|\exp (\alpha ^2n^2-2\alpha (\alpha -1)n) \\
\ge a_1(n+1)^{-b(n+1)-1-\varepsilon }\exp (\alpha ^2n^2-2\alpha (\alpha -1)n)\\
=a_1\exp \left\{ \alpha ^2n^2-2\alpha (\alpha -1)n-(bn+b+1+\varepsilon )\log (n+1)\right\} \to \infty ,
\end{multline*}
as $n\to \infty$. Therefore the series (43) does not converge in $\left( \Phi_+^{1-\alpha }\right)'$. $\qquad \blacksquare$

\end{document}